\documentclass[psamsfonts]{amsart}
\usepackage{amssymb,eucal,latexsym,xy,graphics}
\xyoption{all}
\xyoption{ps}

\theoremstyle{plain}

\newtheorem{lemma}{Lemma}[section]
\newtheorem{theorem}[lemma]{Theorem}

\newtheorem{corollary}[lemma]{Corollary}

\newtheorem*{stat}{\name}
\newcommand{\name}{testing}

\theoremstyle{definition}

\newtheorem{problem}{Problem}

\theoremstyle{remark}
\newtheorem{remark}[lemma]{Remark}

\newcommand{\qedc}{{\qed}~{\rm Claim~{\theclaim}.}}

\numberwithin{equation}{section}

\newcommand{\pup}[1]{\textup{(}{#1}\textup{)}}

\newcommand{\xe}{{\boldsymbol{e}}}
\newcommand{\xf}{{\boldsymbol{f}}}

\newcommand{\xuf}{{\mathbf{f}}}

\newcommand{\xu}{{\boldsymbol{u}}}

\newcommand{\xd}{{\mathbf{d}}}
\newcommand{\xt}{\mathbf{t}}

\newcommand{\Alg}{\mathbf{Alg}(\Sigma)}
\newcommand{\ALatt}{\mathbf{ALatt}}

\newcommand{\Bow}{\mathcal{D}_{\scriptscriptstyle{\bowtie}}}

\newcommand{\Pow}{\mathfrak{P}}

\newcommand{\msd}{meet-sem\-i\-dis\-trib\-u\-tive}

\newcommand{\eps}{\varepsilon}
\newcommand{\es}{\varnothing}
\newcommand{\into}{\hookrightarrow}

\newcommand{\seq}[1]{\langle{#1}\rangle}
\newcommand{\famm}[2]{\left\langle#1\mid#2\right\rangle}

\newcommand{\set}[1]{\{#1\}}

\newcommand{\zero}{\mathbf{0}}
\newcommand{\one}{\mathbf{1}}
\newcommand{\two}{\mathbf{2}}

\newcommand{\CC}{\mathbf{C}}
\newcommand{\DD}{\mathbf{D}}

\newcommand{\VV}{\mathbf{V}}

\newcommand{\EE}{\mathcal{E}}

\DeclareMathOperator{\Con}{Con}

\DeclareMathOperator{\im}{im}
\newcommand{\id}{\mathrm{id}}
\newcommand{\jz}{$\langle\vee,0\rangle$}
\newcommand{\jzu}{$\langle\vee,0,1\rangle$}
\newcommand{\jzs}{\jz-semi\-lat\-tice}
\newcommand{\jzus}{\jzu-semi\-lat\-tice}

\newcommand{\mh}{meet-ho\-mo\-mor\-phism}
\newcommand{\jzuh}{\jzu-ho\-mo\-mor\-phism}
\newcommand{\jze}{\jz-em\-bed\-ding}
\newcommand{\jzue}{\jzu-em\-bed\-ding}

\begin{document}

\title[Diagrams of finite Boolean semilattices]{Congruence lifting of
diagrams of finite Boolean semilattices requires large
congruence varieties}

  \author[J.~T\r uma]{Ji\v r\'\i\ T\r uma}
 \address{Charles University in Prague\\
          Faculty of Mathematics and Physics\\
          Department of Algebra\\
          Sokolovsk\'a 83\\
          Charles University\\
          186 00 Praha 8\\
          Czech Republic}
  \email{tuma@karlin.mff.cuni.cz}

\author[F.~Wehrung]{Friedrich Wehrung}
\address{LMNO, CNRS UMR 6139\\ D\'epartement de
Math\'ematiques\\ Universit\'e de Caen\\ 14032 Caen Cedex\\
France}
\email{wehrung@math.unicaen.fr}

\urladdr{http://www.math.unicaen.fr/\~{}wehrung}

\date{\today}



\thanks{Both authors were partially supported by the institutional
grant CEZ:J13/98:1132000007a and INTAS project 03-51-4110. The first
author was partially supported by grants GA UK 284/2003 and GA CR
201/02/0594, the second author by FRVS 2731/2003 and by the Fund of
Mobility of the Charles University (Prague)}

\begin{abstract}
We construct a diagram $\Bow$, indexed by a finite partially
ordered set, of finite Boolean \jzus s and \jzue s, with top
semilattice $\two^4$, such that for any variety~$\VV$ of
algebras, if $\Bow$ has a lifting, with respect to the congruence
lattice functor, by algebras and homomorphisms in~$\VV$, then there
exists an algebra $U$ in~$\VV$ such that the congruence lattice of
$U$ contains, as a $0$,$1$-sublattice, the five-element modular
nondistributive lattice $M_3$. In particular, $\VV$ has an algebra
whose congruence lattice is neither join- nor \msd. Using earlier work
of K.\,A. Kearnes and \'A.~Szendrei, we also deduce that~$\VV$ has no
nontrivial congruence lattice identity.

In particular, there is no functor $\Phi$ from finite Boolean
semilattices and \jzue s to lattices and lattice embeddings such that
the composition $\Con\Phi$ is equivalent to the identity (where
$\Con$ denotes the congruence lattice functor), thus
solving negatively a problem raised by P. Pudl\'ak in 1985 about
the existence of a functorial solution of the Congruence Lattice
Problem.
\end{abstract}

\maketitle

\section{Introduction}\label{S:Intro}

The Congruence Lattice Problem, CLP in short, asks whether every
distributive algebraic lattice is isomorphic to the congruence lattice
of a lattice. Most of the recent efforts aimed at solving this problem
are focusing on lifting not only individual semilattices, but
\emph{diagrams} of semilattices, with respect to the functor $\Con$,
that with a lattice associates its congruence lattice. This approach
has been initiated by P. Pudl\'ak in \cite{Pudl}. For more
information, we refer the reader to the survey papers by G. Gr\"atzer
and E.\,T. Schmidt \cite{GrScC} and by the authors of the present
paper \cite{TuWe2}.

As observed by the authors in \cite{TuWe1}, the diagram of
Figure~\ref{Fig:NoHomLift}, labeled by Boolean semilattices and
\jzuh s,
\begin{figure}[htb]
 \[
 {
 \def\labelstyle{\displaystyle}
 \xymatrix{
 \two^2\ar@<.5ex>[r]^{\pi} &
 \two\ar@<.5ex>[l]^{\eps}\ar@(ur,dr)[]|{\id}
 }}
 \]
\caption{A diagram unliftable by any algebras.}
\label{Fig:NoHomLift}
\end{figure}
where $\pi(x,y)=x\vee y$ and $\eps(x)=\seq{x,x}$, cannot be lifted
by lattices, see \cite[Theorem~8.1]{TuWe1}. In fact, the proof of
\cite[Theorem~8.1]{TuWe1} shows that there is no lifting of this
diagram by \emph{algebras}, and so this is a result of
universal algebra! However, the map
$\pi$ is not one-to-one, and so the problem remained open whether any
diagram of finite Boolean semilattices and \emph{\jze s} could be
lifted. This problem was first raised by P. Pudl\'ak in 1985, see
the bottom of Page~100 in \cite{Pudl}. It was later attacked by the
authors of the present paper, see Problems~1 and~2 in
\cite{TuWe1}. We shall refer to this problem as \emph{Pudl\'ak's
problem}.

In general, it is not known whether it is decidable,
for a given finite diagram, to have a lifting by, say, lattices, and
tackling even quite simple diagrams may amount to considerable work
with \emph{ad hoc} methods greatly varying from one diagram to the other;
see, for example, the cube diagrams (indexed by $\two^3$) constructed in
\cite{TuWe1}. This partly explains why Pudl\'ak's problem
had remained open for such a long time.

In the present paper, we solve Pudl\'ak's problem by
the negative, by constructing a diagram, indexed by a finite poset,
of finite Boolean semilattices and \jzue s, denoted by $\Bow$, that
cannot be lifted, with respect to the $\Con$ functor, by lattices,
see Theorem~\ref{T:NonLiftM3}. Again, it turns out that as for the
diagram of Figure~\ref{Fig:NoHomLift}, lattice structure does not
really matter, and our negative result holds in any variety satisfying
a non-trivial congruence lattice identity. Whether it holds in any
variety of algebras is still an open problem. 

As the top semilattice of $\Bow$ is $\Pow(4)$, this makes it the
``shortest'' (in terms of the top semilattice) diagram of finite Boolean
\jzus s and \jzue s that cannot be lifted by lattices.

\section{Basic concepts}\label{S:NotTerm}

We denote by $\ALatt$ the category of algebraic lattices and
compactness-pre\-serv\-ing complete join-homomorphisms, and by $\Alg$
the category of algebras of a given similarity type~$\Sigma$ with
$\Sigma$-homomorphisms. For algebras $A$ and $B$ and a homomorphism
$f\colon A\to B$, we denote by $\Con f$ the map that with every
congruence~$\alpha$ of $A$ associates the congruence of $B$
generated by all pairs
$\seq{f(x),f(y)}$, where $\seq{x,y}\in\alpha$. Then $\Con f$ is a
compactness-preserving complete join-homomorphism. The
correspondence $A\mapsto\Con A$,
$f\mapsto\Con f$ is a functor from $\Alg$ to $\ALatt$, that we shall
denote by $\Con$ and call the \emph{congruence lattice functor} on
$\Alg$, see \cite[Section~5.1]{TuWe1}. For an algebra
$A$, we denote by $\zero_A$ the identity congruence of $A$ and by
$\one_A$ the coarse congruence of $A$.

A \emph{diagram} of a category $\DD$, indexed by a category $\CC$, is a
functor from $\CC$ to~$\DD$. Most of our diagrams will be
indexed by posets, the latter being viewed as categories in which
hom sets have at most one element. A \emph{lifting} of a diagram
$\Phi\colon\CC\to\ALatt$ of algebraic lattices is a diagram
$\Psi\colon\CC\to\Alg$ of algebras such that the functor~$\Phi$ and
the composition $\Con\Psi$ are naturally equivalent.

We put $n=\set{0,1,\dots,n-1}$, for every natural number $n$. For a
set $X$, we denote by $\Pow(X)$ the powerset algebra of $X$.

\section{Varieties satisfying a nontrivial congruence lattice identity}
\label{S:CongId}

We first recall some basic definitions and facts of commutator theory,
see R.\,N. McKenzie, G.\,F. McNulty, and W.\,F.
Taylor \cite[Section~4.13]{MMTa87}. For a congruence $\theta$ of an
algebra $A$ and strings $\vec a=\famm{a_i}{i<n}$ and
$\vec b=\famm{b_i}{i<n}$ of elements of $A$, let
$\vec a\equiv_{\theta}\vec b$ hold, if
$a_i\equiv_{\theta}b_i$ holds for all $i<n$.
For congruences $\alpha$, $\beta$,
and $\delta$ of an algebra $A$, we say that $\alpha$
\emph{centralizes $\beta$ modulo $\delta$}, in symbol
$C(\alpha,\beta;\delta)$, if for all $m$, $n<\omega$ and every
$(m+n)$-ary term $t$ of the similarity type of $A$,
 \[
 t(\vec a,\vec p)\equiv_{\delta}t(\vec a,\vec q)\ \Longrightarrow\ 
 t(\vec b,\vec p)\equiv_{\delta}t(\vec b,\vec q),
 \]
for all $\vec a$, $\vec b\in A^m$ and $\vec p$, $\vec q\in A^n$ such
that $\vec a\equiv_{\alpha}\vec b$ and $\vec p\equiv_{\beta}\vec q$.

We state in the following lemma a few standard properties of the
relation $C(\alpha,\beta;\delta)$, see \cite[Lemma~4.149]{MMTa87}.

\begin{lemma}\label{L:BasicComm}
Let $A$ be an algebra.
\begin{enumerate}
\item For all $\alpha$, $\beta\in\Con A$, there exists a
least $\delta\in\Con A$ such that $C(\alpha,\beta;\delta)$. We denote
this congruence by $[\alpha,\beta]$, the \emph{commutator} of $\alpha$
and $\beta$.

\item For all $\beta$, $\delta\in\Con A$, there exists a
largest $\alpha\in\Con A$ such that
$C(\alpha,\beta;\delta)$. We denote this congruence by $(\delta:\beta)$,
the \emph{centralizer} of $\delta$ modulo $\beta$.

\item $[\alpha,\beta]\subseteq\alpha\cap\beta$, for all $\alpha$,
$\beta\in\Con A$.
\end{enumerate}
\end{lemma}

We say that an algebra $A$ is \emph{Abelian}, if
$[\one_A,\one_A]=\zero_A$; equivalently,
$(\zero_A:\nobreak\one_A)=\nobreak\one_A$. A \emph{weak difference term
\pup{resp. weak difference polynomial}} on
$A$ (see K.\,A. Kearnes and \'A. Szendrei \cite{KeSz98}) is a ternary
term (resp., polynomial) $\xd$ such that
 \[
 \xd(x,y,y)\equiv_{[\theta,\theta]}x\equiv_{[\theta,\theta]}\xd(y,y,x),
 \]
for all $\theta\in\Con A$ and all $\seq{x,y}\in\theta$. A weak
difference term for a variety $\VV$ is a ternary term that is a weak
difference term in any algebra of $\VV$.
The statement of the following lemma has been pointed to the authors by
Keith Kearnes.

\begin{lemma}\label{L:WDterm}
Let $A$ be an algebra with a weak difference polynomial $\xd(x,y,z)$. If
there exists a $0$,$1$-homomorphism of $M_3$ to $\Con A$, then $A$ is
Abelian.
\end{lemma}

\begin{proof}
By assumption, there are $\alpha$, $\beta$, $\gamma\in\Con A$ such that
$\alpha\cap\beta=\alpha\cap\gamma=\beta\cap\gamma=\zero_A$ and
$\alpha\vee\beta=\alpha\vee\gamma=\beta\vee\gamma=\one_A$. It follows
from Lemma~\ref{L:BasicComm}(iii) that
$[\beta,\alpha]=[\gamma,\alpha]=\zero_A$, which can be written
$\one_A=\beta\vee\gamma\leq(\zero_A:\alpha)$ (see
Lemma~\ref{L:BasicComm}(ii)), so that $[\one_A,\alpha]=\zero_A$;
\emph{a fortiori}, $[\alpha,\alpha]=\zero_A$. Similarly,
$[\beta,\beta]=[\gamma,\gamma]=\zero_A$.

Now let $\seq{x,y}\in\alpha\circ\beta$. Pick
$z\in A$ such that $x\equiv_{\alpha}z\equiv_{\beta}y$. Hence,
 \begin{align*}
 x&=\xd(x,z,z)&&(\text{because }x\equiv_{\alpha}z\text{ and }
 [\alpha,\alpha]=\zero_A)\\
 &\equiv_{\beta}\xd(x,z,y)&&(\text{because }y\equiv_{\beta}z\text{ and }
 d\text{ is a polynomial})\\
 &\equiv_{\alpha}\xd(z,z,y)
 &&(\text{because }z\equiv_{\alpha}x\text{ and }
 d\text{ is a polynomial})\\
 &=y&&(\text{because }y\equiv_{\beta}z\text{ and }
 [\beta,\beta]=\zero_A),
 \end{align*}
and hence $\seq{x,y}\in\beta\circ\alpha$. Therefore, by symmetry, we
have proved that the congruences $\alpha$, $\beta$, and $\gamma$ are
pairwise permutable. The conclusion follows from
\cite[Lemma~4.153]{MMTa87}.
\end{proof}

We say that a variety $\VV$ has \emph{no nontrivial congruence lattice
identity}, if any lattice identity that holds in the congruence lattices
of all algebras in $\VV$ holds in all lattices.

The following lemma sums up a few deep results of universal
algebra, namely Corollaries~4.11 and~4.12 in \cite{KeSz98}. Recall that
an algebra $A$ is \emph{affine}, if there are a ternary term
operation $\xt$ of $A$ and an Abelian group operation
$\seq{x,y}\mapsto x-y$ on~$A$ such that $\xt(x,y,z)=x-y+z$ for all $x$,
$y$, $z\in A$ and $\xt$ is a homomorphism from~$A^3$ to~$A$.
For further information we refer the reader to C. Herrmann \cite{Herr79}
or R. Freese and R.\,N. McKenzie \cite[Chapter~5]{Comm}. We call $\xt$ a
\emph{difference operation} for $A$.

\goodbreak

\begin{theorem}[K.\,A. Kearnes and \'A. Szendrei]\label{T:StrEq}
For a variety $\VV$ satisfying a nontrivial congruence lattice identity,
the following statements hold.
\begin{enumerate}
\item $\VV$ has a weak difference term.

\item Any Abelian algebra of $\VV$ is affine \pup{thus it has
permutable congruences}.
\end{enumerate}
\end{theorem}

In fact, it follows from \cite[Theorem~4.8]{KeSz98} that (i) implies
(ii) in Theorem~\ref{T:StrEq}. This observation is used in
Remark~\ref{Rk:doingless}.

An algebra $A$ is \emph{Hamiltonian}, if every subalgebra of $A$ is a
congruence class of~$A$.

\begin{lemma}[folklore]\label{L:Aff2Ham}
Every affine algebra $A$ is Hamiltonian.
\end{lemma}

\begin{proof}
Let $\xt$ be a difference operation for $A$, with an Abelian
group operation $-$ on~$A$, and let $U$ be a subalgebra
of $A$. For all $x\in A$ and all $u$, $v\in U$, it follows from the
equation $x+v=\xt(x+u,u,v)$ that $x+u\in U$ implies that $x+v\in U$.
Hence, picking $u\in U$, we can define a binary relation $\equiv$ on $A$
by setting
 \begin{equation}\label{Eq:txyequiv}
 x\equiv y\ \Longleftrightarrow\ \xt(x,y,u)\in U,
 \quad\text{for all }x,\,y\in A,
 \end{equation}
and the relation $\equiv$ is independent of the choice of $u$.
The relation $\equiv$ is obviously reflexive. Now observe that the
equation $\xt(y,x,u)=\xt(u,\xt(x,y,u),u)$ holds for all $x$, $y\in A$. In
particular, if $x\equiv y$, then, as $U$ is closed under $\xt$ (because
$\xt$ is a term operation of $A$), $\xt(y,x,u)$ belongs to $U$, that is,
$y\equiv x$. Hence $\equiv$ is symmetric.
Similarly, by using the equation
$\xt(x,z,u)=\xt(\xt(x,y,u),u,\xt(y,z,u))$, that holds for all $x$, $y$,
$z\in A$, we obtain that $\equiv$ is transitive. By using the assumption
that~$\xt$ is a homomorphism from
$A^3$ to $A$ and by the independence of \eqref{Eq:txyequiv} from $u$,
we obtain that $\equiv$ is compatible with all the operations of $A$;
whence it is a congruence of~$A$. Finally, it is obvious that $U$ is
the equivalence class of $u$ modulo $\equiv$.
\end{proof}

We shall use Lemma~\ref{L:Aff2Ham} in its following special form:
\emph{every subalgebra of a simple affine algebra is either a
one-element algebra or the full algebra}.

\section{The diagram $\Bow$}\label{S:FunctSol4}

We define \jzuh s $\xe\colon\Pow(1)\to\Pow(2)$,
$\xf_i\colon\Pow(2)\to\Pow(3)$, and
$\xu_i\colon\Pow(3)\to\Pow(4)$ (for $i<3$) by their values on atoms:
 \begin{gather*}
 \xe\colon\set{0}\mapsto\set{0,1};\\
 \xf_0\colon
 \begin{cases}
 \set{0}&\mapsto\set{0,1}\\
 \set{1}&\mapsto\set{0,2}
 \end{cases},\qquad
 \xf_1\colon
 \begin{cases}
 \set{0}&\mapsto\set{0,1}\\
 \set{1}&\mapsto\set{1,2}
 \end{cases},\qquad
 \xf_2\colon
 \begin{cases}
 \set{0}&\mapsto\set{0,2}\\
 \set{1}&\mapsto\set{1,2}
 \end{cases},\\
 \xu_0\colon\begin{cases}
 \set{0}&\mapsto\set{0}\\
 \set{1}&\mapsto\set{1,3}\\
 \set{2}&\mapsto\set{2,3}\\
 \end{cases},\qquad
 \xu_1\colon\begin{cases}
 \set{0}&\mapsto\set{0,3}\\
 \set{1}&\mapsto\set{1}\\
 \set{2}&\mapsto\set{2,3}\\
 \end{cases},\qquad
 \xu_2\colon\begin{cases}
 \set{0}&\mapsto\set{0,3}\\
 \set{1}&\mapsto\set{1,3}\\
 \set{2}&\mapsto\set{2}\\
 \end{cases}.
 \end{gather*}
Our diagram, that we shall denote by $\Bow$, is represented on
Figure~\ref{Fig:BowTie4}. All its arrows are \jzue s, and none of them
except $\xe$ is a \mh. It is easy to verify that the diagram $\Bow$
is commutative.
\begin{figure}[htb]
 \[
 {
 \def\labelstyle{\displaystyle}
 \xymatrixrowsep{2pc}\xymatrixcolsep{1.5pc}
 \xymatrix{ & \Pow(4) & \\
 \Pow(3)\ar[ru]^{\xu_0} & \Pow(3)\ar[u]|-(.45){\strut\xu_1} &
 \Pow(3)\ar[lu]_{\xu_2}\\
 & &\\
 \Pow(2)\ar[uu]^{\xf_0}\ar[ruu]^(.35){\xf_0}\ar[rruu]|-(.15){\xf_0}&
 \Pow(2)\ar[luu]|-(.7){\strut\xf_1}\ar[uu]|-(.7){\strut\xf_1}
 \ar[ruu]|-(.7){\strut\xf_1}
 &
 \Pow(2)\ar[lluu]|-(.15){\xf_2}\ar[luu]_(.35){\xf_2}\ar[uu]_{\xf_2}\\
 &\Pow(1)\ar[lu]^{\xe}\ar[u]|-(.45){\strut\xe}\ar[ru]_{\xe} &
 }}
 \]
\caption{The diagram $\Bow$.}\label{Fig:BowTie4}
\end{figure}

\begin{figure}[hbt]
 \includegraphics{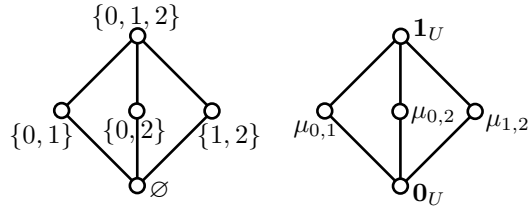}
 \caption{Representing $M_3$ by either subsets or congruences.}
 \label{Fig:M3}
 \end{figure}

We denote as usual by $M_3$ the unique modular nondistributive lattice
with five elements. We shall represent it as a semilattice of subsets
of $3$, by
 \begin{equation}\label{Eq:ReprM3}
 M_3=\set{\es,\set{0,1},\set{0,2},\set{1,2},\set{0,1,2}},
 \end{equation}
see the left hand side of Figure~\ref{Fig:M3}. We are now ready to
formulate our main result.

\begin{theorem}\label{T:NonLiftM3}
Let~$\VV$ be a variety of algebras. If the diagram $\Bow$ admits a
lifting, with respect to the congruence lattice functor, by
algebras and homomorphisms in~$\VV$, then there exists an
algebra~$U$ in~$\VV$ such that $M_3$ has a $0$,$1$-lattice embedding
into $\Con U$. Furthermore, $\VV$ satisfies no nontrivial congruence
lattice identity.
\end{theorem}

The conclusion of Theorem~\ref{T:NonLiftM3} can be further strengthened,
see Remark~\ref{Rk:doingless}.

In particular, Theorem~\ref{T:NonLiftM3} implies immediately that
$\Bow$ cannot be lifted by \emph{lattices}.

\goodbreak

\begin{lemma}\label{L:imuicap}
The semilattice $\im\xu_0\cap\im\xu_1\cap\im\xu_2$ is isomorphic to
$M_3$.
\end{lemma}

\begin{proof}
It is straightforward to verify that for all $X_0$, $X_1$,
$X_2\subseteq 3$, $\xu_0(X_0)=\xu_1(X_1)=\xu_2(X_2)$ if{f}
$X_0=X_1=X_2$ belongs to $M_3$, see \eqref{Eq:ReprM3} and the left hand
side of Figure~\ref{Fig:M3}.
\end{proof}

\section{Proof of Theorem~\ref{T:NonLiftM3}}\label{S:PfThmM3}

\begin{figure}[htb]
 \[
 {
 \def\labelstyle{\displaystyle}
 \xymatrixrowsep{2pc}\xymatrixcolsep{1.5pc}
 \xymatrix{ & B & & & & B &\\
 B_0\ar[ru] & B_1\ar[u] & B_2\ar[lu] & &
 B_0\ar[ru] & B_1\ar[u] & B_2\ar[lu]\\
 & & & & & U\ar[lu]\ar[u]\ar[ru] &\\
 A_0\ar[uu]\ar[ruu]\ar[rruu] & A_1\ar[luu]\ar[uu]\ar[ruu] &
 A_2\ar[lluu]\ar[luu]\ar[uu] & &
 A_0\ar[ru] & A_1\ar[u] & A_2\ar[lu]\\
 & A\ar[lu]\ar[u]\ar[ru] & & & & A\ar[lu]\ar[u]\ar[ru] &
 }}
 \]
\caption{The diagram $\EE$ and the algebra $U$.}\label{Fig:EE4}
\end{figure}
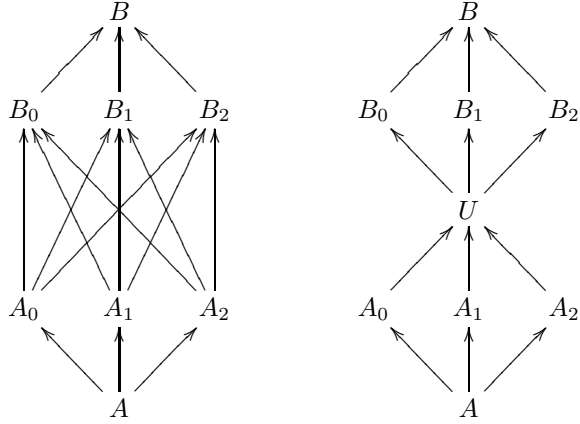

We let a diagram $\EE$ of
algebras in~$\VV$ lift $\Bow$, and we label $\EE$ as on the left hand
side of Figure~\ref{Fig:EE4}.
For vertices~$X$ and~$Y$ of $\EE$, denote by $f_{X,Y}$ the unique
homomorphism from~$X$ to~$Y$ arising from $\EE$. Furthermore,
put $\xuf_{X,Y}=\Con f_{X,Y}$, a \jzuh\ from $\Con X$ to $\Con Y$.
For example, $f_{A_0,B}$ is a homomorphism from $A_0$ to $B$, with
$f_{A_0,B}=f_{B_0,B}\circ f_{A_0,B_0}=f_{B_1,B}\circ f_{A_0,B_1}$, and
so on.

We let the equivalence between $\Con\EE$
and $\Bow$ be witnessed by isomorphisms~$\eps_X$ from the
corresponding vertex of~$\Bow$ onto $\Con X$, for $X$ among $A$,
$A_0$, $A_1$, $A_2$, $B_0$, $B_1$, $B_2$, $B$. For example,
$\eps_{A_i}\colon\two^2\to\Con A_i$,
$\eps_{B_j}\colon\two^3\to\Con B_j$, and
$\eps_{B_j}\circ\xf_i=\xuf_{A_i,B_j}\circ\eps_{A_i}$, for all $i$,
$j<3$.

Since $\EE$ lifts $\Bow$ and all arrows of $\Bow$ are \jze s, so
that, in particular, they separate $0$, all arrows of $\EE$ are
\emph{embeddings}. In particular, we may replace each vertex $X$ of
$\EE$ by its image $\im f_{X,B}$ in $B$, and thus assume that all
$f_{X,Y}$-s in~$\EE$ are (set-theoretical) \emph{inclusion
mappings}. We denote by $U$ the subalgebra of $B$ generated by
$A_0\cup A_1\cup A_2$. Of course, $U$ is contained in
$B_0\cap B_1\cap B_2$, see the right hand side of
Figure~\ref{Fig:EE4}.

\begin{lemma}\label{L:OnePres4}
The equalities $\xuf_{A,U}(\one_A)=\xuf_{A_i,U}(\one_{A_i})=\one_U$
hold, for all $i<3$.
\end{lemma}

\begin{proof}
It suffices to prove that the congruence $\theta=\xuf_{A,U}(\one_A)$ is
equal to $\one_U$. Fix
$a\in\nobreak A$, and put $V=[a]_{\theta}$, the $\theta$-block of
$a$ in $U$. As any pair of elements of $A$ is $\theta$-congruent, $V$ is a
subalgebra of $U$. As $\EE$ lifts $\Bow$ and $\xe$ is unit-preserving,
all maps $\xuf_{A,A_i}$, for $i<3$, are unit-preserving. Hence, for any
$x\in A_0$, we obtain the relation $x\equiv a\pmod{\xuf_{A,A_0}(\one_A)}$,
thus, \emph{a fortiori}, $x\equiv a\pmod{\xuf_{A,U}(\one_A)}$;
whence~$A_0$ is contained in~$V$. Similarly, both $A_1$ and $A_2$ are
contained in $V$. Since $U$ is the subalgebra of $B$ generated by
$A_0\cup A_1\cup A_2$, we obtain $V=U$, that is, $\theta=\one_U$.
\end{proof}

Now we put $\xi_{i,j}=\xuf_{A_i,U}\circ\eps_{A_i}(\set{j})$, for all
$i<3$ and $j<2$. So $\xi_{i,j}$ is a compact congruence of $U$.
Furthermore, it follows from Lemma~\ref{L:OnePres4} that
 \begin{equation}\label{Eq:xijoin21}
 \xi_{i,0}\vee\xi_{i,1}=\one_U,\text{ for all }i<3.
 \end{equation}
Now we put
 \begin{equation}\label{Eq:themus}
 \mu_{0,1}=\xi_{0,0}\vee\xi_{1,0},\quad
 \mu_{0,2}=\xi_{0,1}\vee\xi_{2,0},\quad
 \mu_{1,2}=\xi_{1,1}\vee\xi_{2,1}.
 \end{equation}

\begin{lemma}\label{L:M3intoConU}
There exists a unique $0$,$1$-lattice embedding from $M_3$ into
$\Con U$ that sends $\set{i,j}$ to $\mu_{i,j}$, for all $i<j<3$.
\end{lemma}

\begin{proof}
It follows from \eqref{Eq:xijoin21} that
$\mu_{0,1}\vee\mu_{0,2}=\mu_{0,1}\vee\mu_{1,2}=\mu_{0,2}\vee\mu_{1,2}
=\one_U$. Now we need to prove that any two $\mu_{i,j}$-s meet at
zero. We observe that for all $i$, $k<3$ and all $j<2$,
 \[
 \xuf_{U,B_k}(\xi_{i,j})=
 \xuf_{U,B_k}\circ\xuf_{A_i,U}\circ\eps_{A_i}(\set{j})=
 \xuf_{A_i,B_k}\circ\eps_{A_i}(\set{j})=
 \eps_{B_k}\circ\xf_i(\set{j}).
 \]
This makes it easy to calculate the values $\xuf_{U,B_k}(\mu_{i,j})$,
for $i<j<3$ and $k<3$. We obtain, using \eqref{Eq:themus},
 \[
 \xuf_{U,B_k}(\mu_{i,j})=\eps_{B_k}(\set{i,j}),
 \text{ for all }i<j<3\text{ and }k<3.
 \]
It follows immediately that
$\xuf_{U,B_k}(\mu_{0,1}\wedge\mu_{0,2})=\eps_{B_k}(X_k)$, for some
$X_k\subseteq\set{0}$. Hence, applying the observation that
$\xuf_{B_k,B}\circ\xuf_{U,B_k}=\xuf_{U,B}$ (independent of $k$), we
obtain that $\xu_0(X_0)=\xu_1(X_1)=\xu_2(X_2)$, hence, by
Lemma~\ref{L:imuicap}, $X_0=X_1=X_2$ is not a singleton. Therefore,
$X_0=\es$, so $\xuf_{U,B_0}(\mu_{0,1}\wedge\mu_{0,2})=\zero_{B_0}$.
Since $\xuf_{U,B_0}$ separates zero, it follows that
$\mu_{0,1}\wedge\mu_{0,2}=\zero_U$. The proofs that
$\mu_{0,1}\wedge\mu_{1,2}=\zero_U$ and
$\mu_{0,2}\wedge\mu_{1,2}=\zero_U$ are similar. Therefore, the right
hand side of Figure~\ref{Fig:M3} represents a $0$,$1$-sublattice of
$\Con U$ isomorphic to~$M_3$.
\end{proof}

As the subalgebra $U$ of $B$ belongs to~$\VV$,
the result of Lemma~\ref{L:M3intoConU} completes the proof of
the first part of Theorem~\ref{T:NonLiftM3}.

Now suppose that $\VV$ (our variety lifting $\Bow$) satisfies a
nontrivial congruence lattice identity. By Theorem~\ref{T:StrEq}(i),
$\VV$ has a weak difference term. By Lemmas~\ref{L:WDterm}
and~\ref{L:M3intoConU}, the algebra $U$ is Abelian, thus so is $A_0$.
By Theorem~\ref{T:StrEq}(ii), $A_0$ is affine. Since
$\Con A_0\cong\two^2$ and $A_0$ has permutable congruences, we obtain (up
to isomorphism) that $A_0=A'\times A''$, for simple algebras~$A'$
and~$A''$, and $\xe$ is lifted by an embedding of the form
$x\mapsto\seq{e'(x),e''(x)}$, for embeddings $e'\colon A\into A'$ and
$e''\colon A\into A''$. By Lemma~\ref{L:Aff2Ham}, both $A'$ and $A''$ are
Hamiltonian, thus both~$e'$ and~$e''$ are isomorphisms, whence
$A'\cong A''$. Now observe that any term giving a difference operation
on the affine algebra $A_0$ satisfies Mal'cev's equations.
Hence, taking $A'=A''$ and using \cite[Lemma~4.154]{MMTa87}, we obtain
that the smallest congruence of~$A_0$ collapsing the diagonal of~$A'$ is
a complement, in $\Con A_0$, of both projection kernels in
$A_0=A'\times A'$, which contradicts
$\Con A_0\cong\two^2$. This concludes the proof of
Theorem~\ref{T:NonLiftM3}.

\begin{remark}\label{Rk:doingless}
As $M_3$ does not embed into the congruence lattice of any lattice, a
mere solution of Pudl\'ak's problem does not require any use of commutator
theory.

On the other hand, it follows from \cite[Theorem~4.8]{KeSz98} that (i)
implies (ii) in the statement of Theorem~\ref{T:StrEq}. Therefore, the
conclusion of Theorem~\ref{T:NonLiftM3} can be strengthened into saying
that $\VV$ has no weak difference term. However, that this is
indeed a strengthening is not trivial, and it follows from the deep
\cite[Corollary~4.12]{KeSz98}.
\end{remark}

\section{Discussion}\label{S:Conseq}

A first immediate corollary of Theorem~\ref{T:NonLiftM3} is the
following.

\begin{corollary}\label{C:NonLift1}
Let~$\VV$ be a variety of algebras with a nontrivial congruence lattice
identity. Then there is no functor $\Phi$ from the category of finite
Boolean \jzus s and \jzue s to~$\VV$ such that the composition
$\Con\Phi$ is equivalent to the identity.
\end{corollary}

The assumption of Corollary~\ref{C:NonLift1} holds, in particular,
for~$\VV$ being the variety of all lattices, which is
congruence-distributive. Therefore, this solves negatively Pudl\'ak's
problem. Of course, the level of generality obtained by
Theorem~\ref{T:NonLiftM3} goes far beyond the failure of
congruence-distributivity---for example, it includes
congruence-modularity.

These results raise the problem whether there is a `quasivariety
version' of the variety result stated in
Corollary~\ref{C:NonLift1}.

\begin{problem}\label{Pb:QVar}
Let~$\VV$ be a variety of algebras. If every finite
poset-indexed diagram of finite Boolean semilattices and \jzue s can be
lifted, with respect to the congruence lattice functor, by a diagram
in~$\VV$, then can every finite lattice be embedded into the
congruence lattice of some algebra in~$\VV$?
\end{problem}

A possibility would be to introduce more complicated variants of the
diagram~$\Bow$, which would yield a sequence $\famm{S_i}{i<\omega}$
of finite lattices, each of which would play a similar role as $M_3$ in
the proof of Lemma~\ref{L:M3intoConU}, and that would generate the
quasivariety of all lattices. However, the main difficulty of the crucial
Lemma~\ref{L:M3intoConU} is the preservation of meets, for which the
lattice $M_3$ is quite special. Although we know how to extend the method
to many finite lattices, we do not know how to get all of them.

Another natural question is whether there is any variety at
all satisfying the natural strengthening of the assumption of
Theorem~\ref{T:NonLiftM3}.

\begin{problem}\label{Pb:ExLarVar}
Does there exist a variety~$\VV$ of algebras such that every finite
poset-indexed diagram of finite Boolean \jzus s and \jzue s can be
lifted, with respect to the congruence lattice functor, by algebras
in~$\VV$?
\end{problem}

Of course, a similar question can be formulated for
finite Boolean \jzs s and \jze s.

In view of Bill Lampe's results \cite{Lamp82}, a natural candidate
for~$\VV$ would be the variety of all \emph{groupoids} (i.e., sets
with a binary operation). By the second author's results in
\cite{Wehr1,Wehr2}, any variety satisfying the conclusion of
Problem~\ref{Pb:ExLarVar} (provided there is any) has the
property that every diagram of finite \emph{distributive} \jzus s
and \jzue s can be lifted, with respect to the congruence lattice
functor, by algebras in~$\VV$. In particular, every distributive
algebraic lattice \emph{with compact unit} would be isomorphic to the
congruence lattice of some algebra in~$\VV$. The latter conclusion
is known (even without the distributivity restriction!) in case~$\VV$
is the variety of all groupoids, see~\cite{Lamp82}. Strangely, whether
the corresponding result holds for arbitrary distributive algebraic
lattices is still an open problem. Also, there exists an algebraic
lattice that is not isomorphic to the congruence lattice of any
groupoid---namely, the subspace lattice of any infinite-dimensional
vector space over any uncountable field, see~\cite{FLT79}.

\section*{Acknowledgments}

This work was partially completed while the second author was
visiting the Charles University (Prague). Excellent conditions
provided by the Department of Algebra are greatly appreciated.

The authors are also grateful to Ralph Freese and Keith Kearnes
for enlightening e-mail discussions about congruence varieties.


\begin{thebibliography}{99}

\bibitem{FLT79}
R. Freese, W.\,A. Lampe, and W. Taylor,
\emph{Congruence lattices of algebras of fixed similarity type. I},
Pacific J. Math. \textbf{82}, 59--68 (1979). 

\bibitem{Comm}
R. Freese and R.\,N. McKenzie,
``Commutator Theory for Congruence Modular Varieties'',
London Math. Soc. Lect. Notes Series \textbf{125},
Cambridge University Press, 1987. 227~p.

\bibitem{GLT2}
G. Gr\"atzer, ``General Lattice Theory. Second edition'', new appendices
by the author with B.\,A. Davey, R. Freese, B. Ganter, M. Greferath, P.
Jipsen, H.\,A. Priestley, H. Rose, E.\,T. Schmidt, S.\,E. Schmidt, F.
Wehrung, and R. Wille. Birkh\"auser Verlag, Basel, 1998. xx+663~p.

\bibitem{Herr79}
C. Herrmann,
\emph{Affine algebras in congruence modular varieties},
Acta Sci. Math. (Szeged) \textbf{41} (1979), 119--125.

\bibitem{KeSz98}
K.\,A. Kearnes and \'A. Szendrei,
\emph{The relationship between two commutators},
Internat. J. Algebra Comput. \textbf{8}, no.~4 (1998), 497--531.

\bibitem{GrScC}
G. Gr{\"a}tzer and E.\,T. Schmidt,
\emph{Congruence Lattices}, Appendix~C in \cite{GLT2}, 519--530.

\bibitem{Lamp82}
W.\,A. Lampe,
\emph{Congruence lattices of algebras of fixed similarity type.
II}, Pacific J. Math. \textbf{103} (1982), 475--508.

\bibitem{MMTa87}
R.\,N. McKenzie, G.\,F. McNulty, and W.\,F. Taylor,
``Algebras, Lattices, Varieties. Volume~I.''
The Wadsworth \& Brooks/Cole Mathematics Series. 
Monterey, California: Wadsworth \& Brooks/Cole Advanced Books \&
Software, 1987. xii+361 p.

\bibitem{Pudl}
P. Pudl\'ak,
\emph{On congruence lattices of lattices}, Algebra Universalis
\textbf{20} (1985), 96--114.

\bibitem{TuWe1}
J. T\r{u}ma and F. Wehrung,
\emph{Simultaneous representations of semilattices by lattices
with permutable congruences}, Internat. J. Algebra Comput.
\textbf{11}, no.~2 (2001), 217--246.

\bibitem{TuWe2}
J. T\r{u}ma and F. Wehrung,
\emph{A survey of recent results on congruence lattices of lattices},
Algebra Universalis \textbf{48}, no.~4 (2002), 439--471.

\bibitem{Wehr1}
F. Wehrung,
\emph{Distributive semilattices as retracts of ultraboolean ones:
functorial inverses without adjunction}, preprint.

\bibitem{Wehr2}
F. Wehrung,
\emph{Lifting retracted diagrams with respect to projectable
functors}, preprint.

\end{thebibliography}
\end{document}